\documentclass[mathematics,article,accept,moreauthors,pdftex]{mdpi}
\usepackage{amsmath} 
\newcommand{\atan}{\operatorname{atan}}
\usepackage{pifont}
\newcommand{\cmark}{\ding{51}}
\newcommand{\xmark}{\ding{55}}
\usepackage{kotex}
\firstpage{1} 
\pubvolume{9}
\issuenum{1}
\articlenumber{1317}
\pubyear{2021}
\copyrightyear{2021}
\externaleditor{Academic Editor: Liliya Demidova} 
\datereceived{13 May 2021} 
\dateaccepted{3 June 2021} 
\datepublished{8 June 2021} 
\hreflink{https://doi.org/10.3390/\linebreak math9121317} 

\Title{Not Another Computer Algebra System: Highlighting \emph{wxMaxima} in Calculus}

\TitleCitation{Not Another Computer Algebra System: Highlighting \emph{wxMaxima} in Calculus}

\Author{Natanael Karjanto $^{1,}$*\orcidA{} and Husty Serviana Husain $^{2}$}

\AuthorNames{N. Karjanto and H. S. Husain}

\AuthorCitation{Karjanto, N.; Husain, H.S.}

\address{%
$^{1}$ \quad Department of Mathematics, University College, Natural Science Campus, Sungkyunkwan University 
Suwon~16419, Republic of Korea \\
$^{2}$ \quad Department of Mathematics Education, Faculty of Mathematics and Natural Science Education\linebreak Indonesia Education University, Bandung 40154, Indonesia; serviana@upi.edu}

\corres{Correspondence: natanael@skku.edu}

\abstract{This article introduces and explains a computer algebra system (CAS) \emph{wxMaxima} for Calculus teaching and learning at the tertiary level. The didactic reasoning behind this approach is the need to implement an element of technology into classrooms to enhance students' understanding of Calculus concepts. For many mathematics educators who have been using CAS, this material is of great interest, particularly for secondary teachers and university instructors who plan to introduce an alternative CAS into their classrooms. By highlighting both the strengths and limitations of the software, we hope that it will stimulate further debate not only among mathematics educators and software users but also also among symbolic computation and software developers.}


\begin{document}

\section{Introduction}

A computer algebra system (CAS) is a program that can solve mathematical problems by rearranging formulas and finding a formula that solves the problem, as opposed to just outputting the numerical value of the result. \emph{Maxima} is a full-featured open-source CAS: the software can serve as a calculator,  provide analytical expressions, and perform symbolic manipulations. Furthermore, it offers a range of numerical analysis methods for equations or systems of equations that otherwise cannot be solved analytically. It can sketch graphical objects with excellent quality. 

What is \emph{wxMaxima}, then? \emph{wxMaxima} is a document-based graphical user interface (GUI) for the CAS \emph{Maxima}. It allows us for using all of \emph{Maxima}'s functions. Additionally, it provides convenient wizards for accessing the most commonly used features, including inline plots and simple animations. Similar to \emph{Maxima}, \emph{wxMaxima} is free of charge, and it is released and distributed under the terms of the GNU General Public License (GPL). This allows for everyone to modify and distribute it, as long as its license remains unmodified. In this article, we use the term ``\emph{wxMaxima}'' more often, but the terms ``\emph{Maxima}'' and ``\emph{wxMaxima}'' can be used interchangeably. 

\emph{Maxima} is different from other well-known so-called \emph{3M} mathematical software (\emph{Maple}, \emph{Matlab}, \emph{Mathematica}), as they are commercial and one needs to purchase a license before using them. Other open-source mathematics software include \emph{Axiom}, \emph{Reduce}, \emph{SageMath}, \emph{Octave} and \emph{Scilab} (both are for numerical computation), \emph{R} (for statistical computing), and \emph{GeoGebra} (for interactive geometry and algebra), where the latter is quite well-known globally among mathematics educators.

Apart from being free and easy to install, \emph{Maxima} is also updated continuously. Currently, \emph{Maxima} can run natively without emulation on the following operating systems: \emph{Windows}, \emph{Mac OS X}, \emph{Linux}, \emph{Berkeley Software Distribution (FreeBSD)}, \emph{Solaris}, and \emph{Android}. An executable file can be downloaded from \emph{Maxima}'s website~\cite{maxima}. In particular, the installation file for \emph{Windows} operating system is available for download at an open-source software community resource \emph{SourceForge}~\cite{sourceforge}. One can simply double-click the executable file and follow the instruction accordingly. After the installation is completed, the software is ready to be launched. The whole process takes less than three minutes in total, depending on the Internet connection speed.

This software is introduced because it is free and under the GPL. As a comparison, the total combined cost for purchasing \emph{3M} software is almost USD~5500 for an educational license (see Table~\ref{cost}). Although the expenditure for personal and student licenses are much lower than those for commercial and professional licenses, for colleagues and practitioners in many developing countries, the price is still considered to be costly. Having this in mind, we promote an open-source software that benefits many people who have limited resources, particularly in less affluent countries. 

For teaching and learning mathematics, \emph{Maxima} is fairly accessible by many people. Although \emph{SageMath} is popular among university professors for teaching Calculus and Linear Algebra thanks to its user-friendly cloud, the server is rather slow, particularly if one attempts to access it from a developing country with modest Internet connectivity. \emph{SageMath} can be downloaded and installed locally, but it is a huge file. Thus, it is another hindrance for many colleagues in developing countries.
\begin{specialtable}[H]
\caption{Software packages and the individuals or institutions through which they were first developed, the year they were launched, and estimated cost in USD.}	\label{cost}
\begin{tabular*}{0.75\textwidth}{@{\extracolsep{\fill}} l l c c @{}}
\toprule
\textbf{Software}	& \textbf{Creator }		 & \textbf{Launched}	& \textbf{Cost (USD)} \\ 
\midrule
\emph{Axiom} 		& Richard Jenks 		 & 1977		& Free \\
\emph{Magma}		& University of Sydney	 & 1990 	& \ USD 1140 \\
\emph{Maple}		& University of Waterloo & 1980		& \ USD 2390 \\
\emph{Mathematica}	& Wolfram Research 		 & 1986 	& \ USD 2495 \\
\emph{Maxima}		& Bill Schelter et al.	 & 1976		& Free \\
\emph{Matlab}		& MathWorks				 & 1989		& \ USD 3150 \\
\emph{SageMath}		& William  A. Stein	 	 & 2005		& Free \\
\bottomrule
\end{tabular*}
\end{specialtable}

In what follows, we cover a literature study on \emph{Maxima} for teaching and learning. A study from Malaysia suggests that students who are exposed to \emph{Maxima} while learning Calculus had a significantly better academic performance as compared with the group that followed a traditional teaching method, and showed a better motivation and more confidence towards the subject~\cite{ayub}. Another example comes from Spain, where Garc\'{i}a et al. proposed to replace \emph{Derive} using \emph{Maxima}~\cite{garcia}. D\'{i}az et al. analyzed the role of \emph{Maxima} in learning Linear Algebra in the context of learning on the basis of competencies~\cite{diaz}. Fedriani and Moyano proposed using \emph{Maxima} in teaching mathematics for business degrees and A-level students~\cite{fedriani}. The authors also presented a report of the main strengths and weaknesses of this software when used in the classroom. Additionally, the CAS \emph{wxMaxima} is used for training future mathematics teachers in Ukraine~\cite{velychko}.

Advanced mathematics can also be explored using \emph{Maxima}, as demonstrated by Dehl~\cite{dehl}. A new possibility for interactive teaching in engineering module using \emph{Maxima} was discussed by \v{Z}\'{a}kov\'{a}~\cite{zakova}. There also exist free Calculus electronic textbooks incorporating \emph{wxMaxima} developed by Zachary Hannan from Solano Community College, Fairfield, California~\cite{hannan}. These books could certainly be adopted into Calculus classrooms. Currently, Hannan is working on PreCalculus, Multivariable Calculus (MVC), Linear Algebra, and Differential Equations textbooks that utilize the software. Thus, we hope to see more \emph{Maxima}-based mathematics textbooks, and this is good news for mathematics educators who are interested in embedding technology and CAS into their classrooms. 

Beyond mathematics, some authors used \emph{Maxima} successfully in Classical Mechanics~\cite{timberlake} and Chemistry~\cite{senese}. In particular, for an efficient path in understanding \emph{Maxima}, the latter provides a thorough introduction to exploiting \emph{Maxima} with the focus on utilizing the \emph{wxMaxima} interface. Woollett provided a series of tutorial notes on \emph{Maxima}. Designed for new users, particularly \emph{Windows} customers, the notes include some nuts-and-bolts suggestions for working with the CAS~\cite{woollett}. Puentedura  designed a \emph{Maxima} tutorial workflow for enhancing and transforming the learning process in science and mathematics. He identified that the CAS has at least three essential roles: as a  number-crunching calculator, as a tool for paper-and-pencil symbolic mathematical derivation, and as a typesetter~\cite{puentedura}.

The list of literature reviews presented above is by no means exhaustive. While the literature offers abundant materials on where \emph{wxMaxima} can outstandingly perform, what it can and cannot do is not entirely clear when it comes to teaching Calculus using the software. This article fills the gap in human--computer interaction, both in the technological and pedagogical senses. Furthermore, by highlighting the software's limitations, we hope to stimulate further debate among the symbolic-computation and mathematics-education communities on how to remedy the situation, perhaps by either providing alternative ways in problem-solving or by improving the technological aspects of the CAS itself. From a pedagogical point of view, the students' feedback that was obtained after implementing \emph{wxMaxima} in teaching and learning for multiple semesters could shed light on its accessibility, digital natives' interaction with technology, and uncover a better way of teaching with technology.

In particular, the limitations of \emph{wxMaxima} are examined by revisiting several examples that were considered in the literature. We present some examples of symbolic integration where \emph{wxMaxima} fails to calculate in most of the cases in relation to results from other software. Interestingly, even other CAS that were often regarded as superior to \emph{wxMaxima} have shortcomings as well. Granted, the field of symbolic integration itself, let alone the area of symbolic computation, particularly with different CAS, is a wide range. Both the cited literature and presented examples are by no means exhaustive lists.

Additionally, although Calculus focuses on the mathematical study of continuous change, the subject in itself contains rich examples of symmetrical objects. For example, even and odd functions are symmetric with respect to the $y$-axis and origin, respectively. Geometrically, the graph of even and odd functions remains unchanged after reflection about the $y$-axis and rotation of 180$^{\circ}$ about the origin, respectively. Consequently, the integral of the former with respect to symmetric intervals is twice the integral from zero to the corresponding upper limit, while the integral of the latter vanishes. We assume that either these functions are integrable and the intervals are finite or the integral converges for infinite intervals. Throughout this article, we consider other examples where symmetry occurs in Calculus and illustrate them using \emph{wxMaxima}.

This article is an extended version of our previous work on embedding technology into Calculus teaching and learning~\cite{karjanto}. We investigate the following research questions: Where can \emph{wxMaxima} perform exceptionally as a CAS? What is its effectiveness? What are \emph{wxMaxima}'s limitations and weaknesses? This paper is organized as follows. After this introduction, Section~\ref{strongweak} discusses \emph{wxMaxima}'s strengths and limitations, which are featured through several examples. Section~\ref{discussion} provides educational benefits when one embeds the software for instructional purposes, in addition to covering some study limitations. Section~\ref{conclusion} \linebreak concludes the study and lays out a future outlook from our discussion.

\section{\emph{wxMaxima}'s Strengths and Weaknesses}		\label{strongweak}

In this section, we cover both \emph{wxMaxima}'s strengths and weaknesses through examples. Although the considered illustrations are mostly Calculus-related, one may consider examples for other subjects as well, including Linear Algebra, Differential Equations, and Discrete Mathematics.

\subsection{\emph{wxMaxima}'s Strengths}

When working with Calculus problems, students can verify the result of their manual computations performed by pencil and paper using \emph{wxMaxima}. Since \emph{wxMaxima} is lightweight and fairly straightforward, simple Calculus calculations can be executed within seconds. These include calculating the limit, finding the derivative of a function, and evaluating both definite and indefinite integrals. The syntax is relatively easy to understand by someone who has no or little experience in programming. In what follows, we consider several \emph{wxMaxima} examples that can be useful for Calculus teaching and learning.	
\begin{Example}	\label{exlim}
Calculating limits, derivatives, and integrals:
\begin{equation*}
\lim\limits_{x \to \,0} \frac{\sin (7x)}{x}, \qquad \qquad
\frac{d}{dx} \cos(3x^2), \qquad \qquad 
\int \frac{1}{1 + x^2} \, dx, \qquad \qquad
\int_0^1 \frac{1}{1 + x^2} \, dx.
\end{equation*}
\end{Example}
Example~\ref{exlim} provides simple examples of calculating limits, derivatives, and integrals, both definite and indefinite. All functions in this example are even. Table~\ref{table1} displays the \emph{wxMaxima} command inputs and their corresponding outputs. Our experience as instructors suggests that \emph{wxMaxima} is also handy for obtaining quick and excellent-quality graphical plots, particularly three-dimensional objects that can be hard to sketch manually. In turn, this visualization enhances students' understanding of applying Calculus concepts.

\begin{specialtable}[H]
\renewcommand{\arraystretch}{1.5}	
\caption{Examples of \emph{wxMaxima} commands and their corresponding outputs covering limit, derivative, integral, and 2D plots.}		\label{table1}
\begin{tabular*}{0.75\textwidth}{@{\extracolsep{\fill}} l l @{}}
\toprule	
\textbf{Input} & \textbf{Output} \\
\midrule	
\verb|(%i1) 'limit(sin(7*x)/x,x,0);|			& ${\displaystyle \lim\limits_{x \to \,0} \frac{\sin (7x)}{x}}$ \\
\verb|(%i2)  limit(sin(7*x)/x,x,0);	|			& $7$ \\
\verb|(%i3) 'diff(cos(3*x^2),x);  |	  	  		& ${\displaystyle \frac{d}{dx} \cos(3x^2)} $ \\
\verb|(%i4)  diff(cos(3*x^2),x);  |	     		& $-6x \sin(3x^2)$ \\
\verb|(%i5) 'integrate(1/(1 + x^2),x); 	|		& ${\displaystyle \int \frac{1}{1 + x^2} \, dx}$ \\
\verb|(%i6) integrate(1/(1 + x^2),x);	|		& $\atan(x)$ \\
\verb|(%i7) integrate(1/(1 + x^2),x,0,1);|		& ${\displaystyle \frac{\pi}{4}}$ \\
\verb|(%i8) plot2d(2*x^3-7*x^2-5*x+4, [x,-2,4.5]); |	& (Figure~\ref{fig2-func}, left panel) \\
\verb|(%i9) plot2d([exp(x), 1 + x], [x,-3,2]); |		& (Figure~\ref{fig2-func}, right panel)\\
\bottomrule
\end{tabular*}
\end{specialtable}

\begin{figure}[H]
\widefigure
\includegraphics[scale=1]{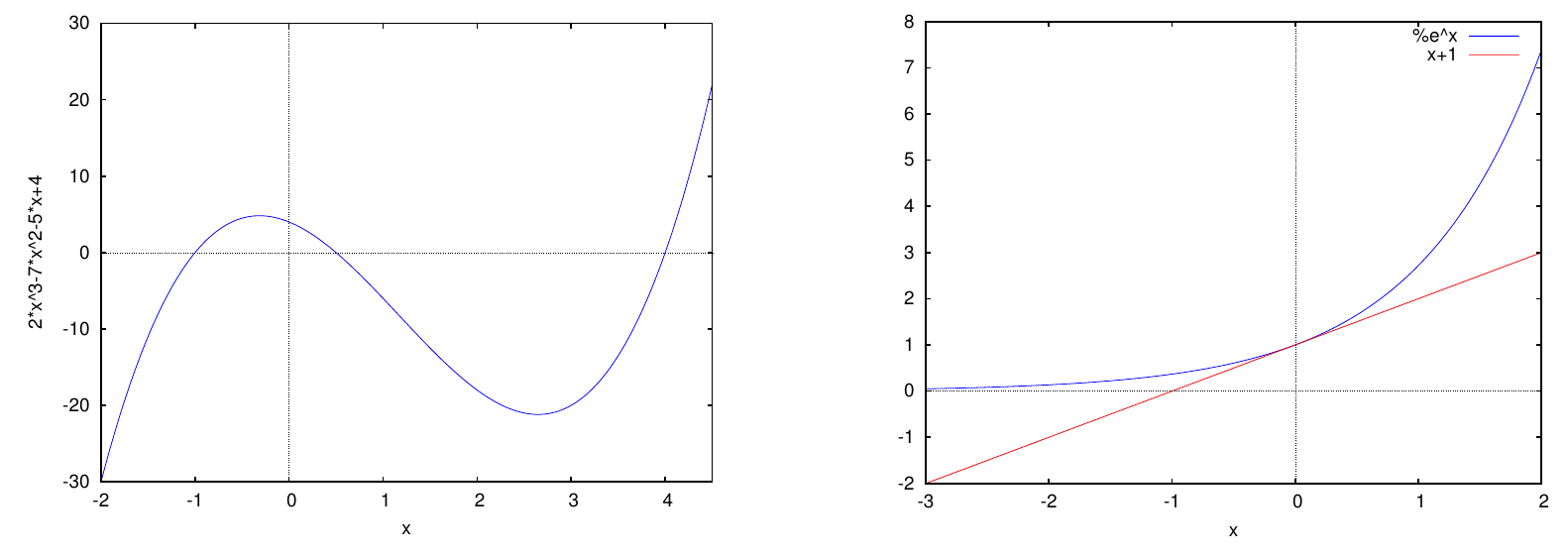} 
\caption{2D plot outputs of a cubic function $y = f(x) = 2x^3 - 7x^2 - 5x + 4$ (\textbf{left}), and the exponential function $y = g(x) = e^x$ and its linear approximation at $(0,1)$, $y = L(x) = 1 + x$ (\textbf{right}).}
\label{fig2-func}
\end{figure}

\begin{figure}[H]
\widefigure
\includegraphics[scale=1]{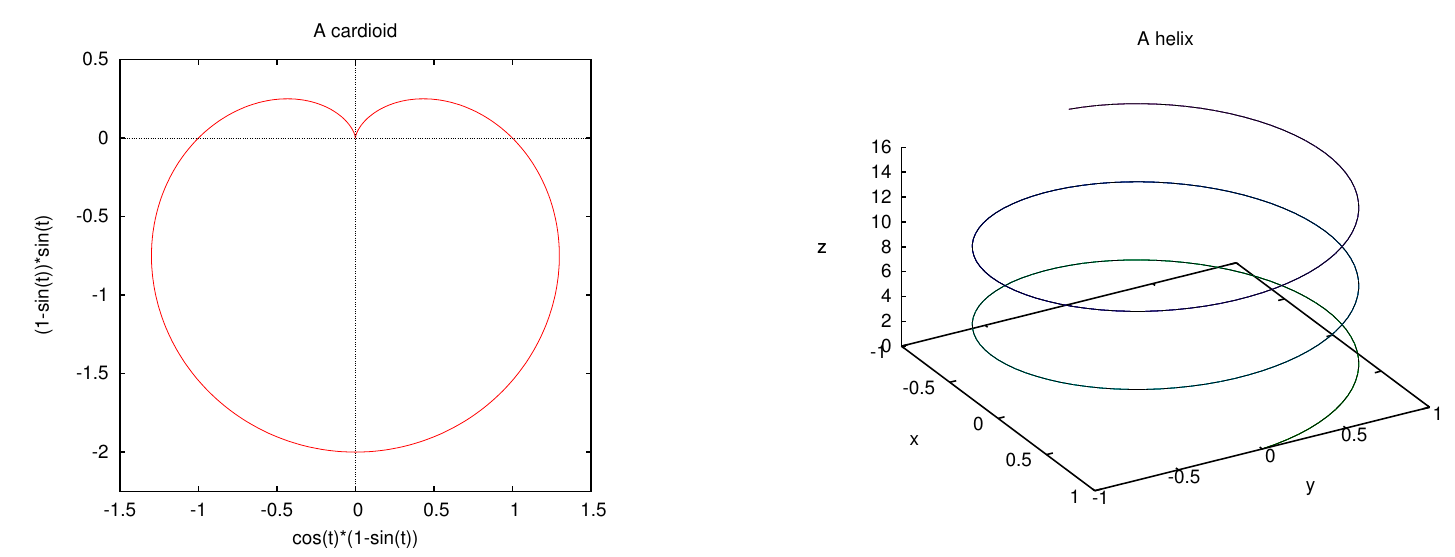}	
\caption{(\textbf{left}) The output of 2D parametric plot showing a polar curve cardioid $r(t) = 1 - \sin t$, $0 \leq t \leq 2\pi$. (\textbf{right}) An output of a 3D parametric plot displaying a vector-valued function helix $\boldsymbol{r}(t) = \langle \cos t, \sin t, t \rangle$, $0 \leq t \leq 5\pi$. These particular cardioid and helix are symmetric with respect to the $y$- and $z$-axes, respectively.} \label{fig3-para}
\end{figure}

\begin{figure}[H]
\widefigure
\includegraphics[scale=1]{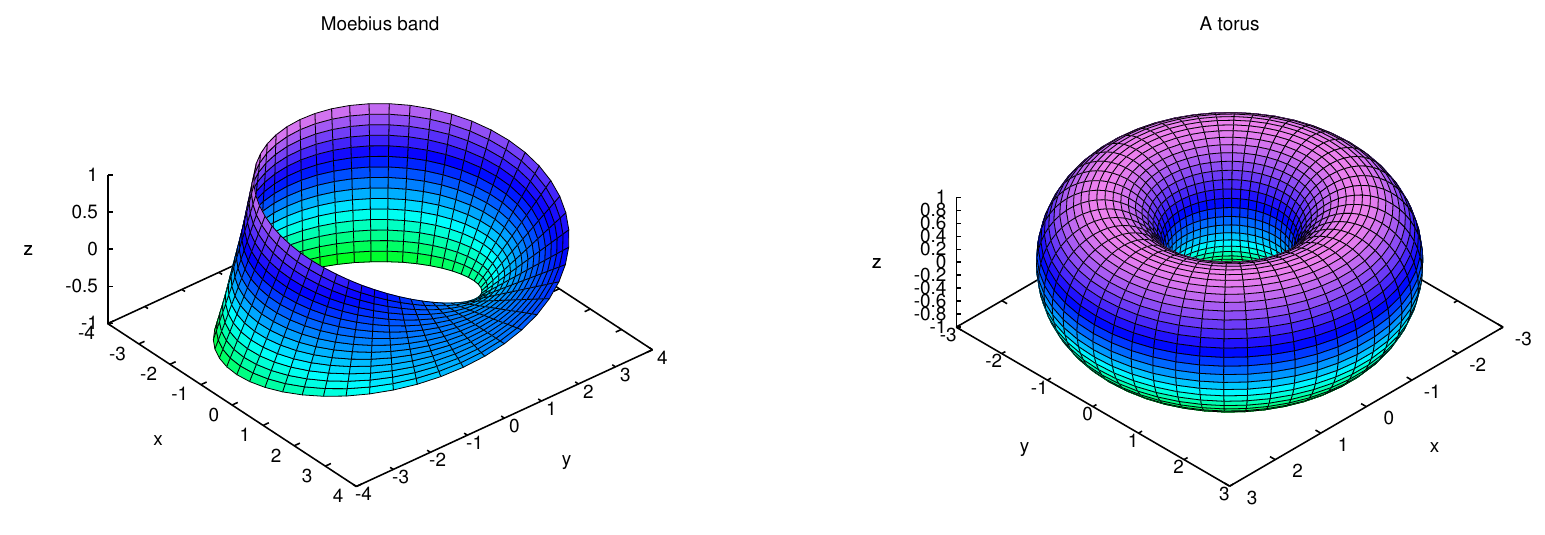}		
\caption{Outputs of (\textbf{left}) the 3D surface M\"obius band $\boldsymbol{s}(x,y) = \left(3 + y \cos \frac{x}{2} \right) \cos x \, \mathbf{i} + \left(3 + y \cos \frac{x}{2}\right) \sin x \, \mathbf{j} + y \sin \frac{x}{2} \, \mathbf{k}$, and (\textbf{right}) the torus $\boldsymbol{t}(\theta,\phi) = (2 + \cos \theta) \cos \phi \, \mathbf{i} + (2 + \cos \theta) \sin \phi \, \mathbf{j} + \sin \theta \, \mathbf{k}$.} \label{fig4-tor}
\end{figure}

\begin{figure}[H]
\widefigure
\includegraphics[scale=1.1]{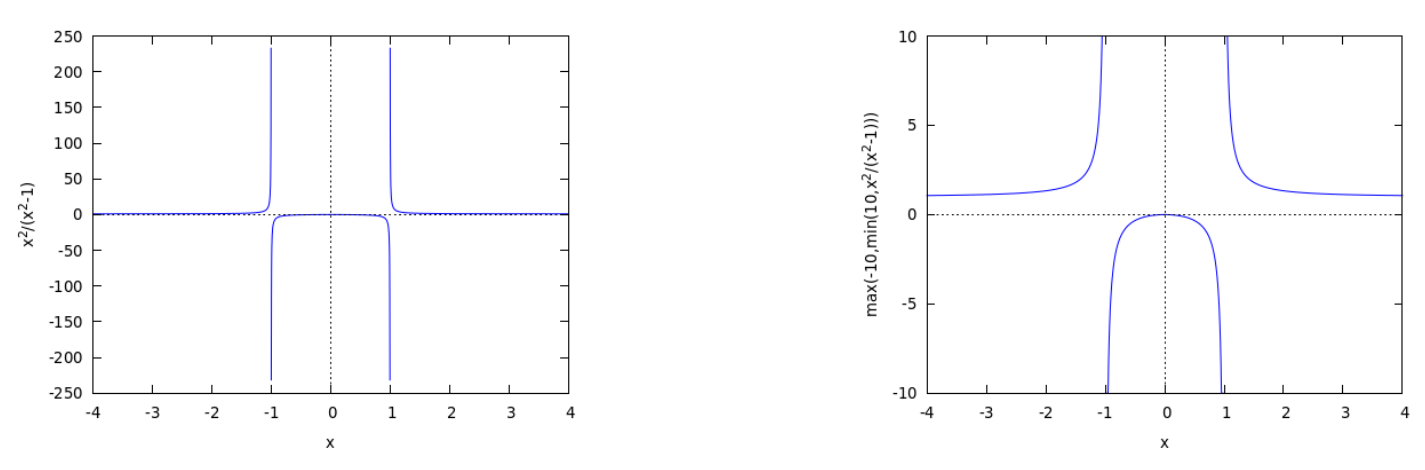}	
\caption{(\textbf{left}) A 2D plot of $y = f(x) = \frac{x^2}{x^2 - 1}$ without any restriction in the vertical scale. (\textbf{right}) The same plot as in the left panel but the vertical axis is restricted to $-10 \leq y \leq 10$.} \label{fig4}
\end{figure}

\begin{Example}		\label{explot}
	Sketching functions (2D plots):
	\begin{align*}
	y = f(x) &= 2x^3 - 7x^2 - 5x + 4  						   & &[\text{cubic function}] \\
	y = g(x) &= e^x \quad \text{and} \quad L(x) = 1 + x \qquad & &[\text{exponential function and its linear approximation}].
	\end{align*}
\end{Example}
Example~\ref{explot} gives examples of generating 2D plots of a single polynomial function and an exponential function together with the corresponding linear approximation. Table~\ref{table1} shows the \emph{wxMaxima} command inputs, and Figure~\ref{fig2-func} displays the corresponding outputs. The plot on the left-hand side can be used to verify the theoretical computation of the local extrema, increasing or decreasing test, and concavity test to a cubic function $f(x) = 2x^3 - 7x^2 - 5x + 4$. The plot on the right-hand side shows a linear approximation $L(x)$ of an exponential function $g(x) = e^x$ at $x = 0$. The sketching quality in \emph{wxMaxima} is remarkably excellent, and the task is completed swiftly. \emph{SageMath} is rather slow to produce plots, and \emph{WolframAlpha}'s free version generates poor quality plots.
\end{paracol}
\nointerlineskip
\unskip
\begin{paracol}{2}
\switchcolumn

\begin{Example}		\label{exparplot}
2D and 3D parametric plots: 	
\begin{align*}
r(t) &= 1 - \sin t, & & \hspace*{-1cm} 0 \leq t \leq 2\pi & &[\text{cardioid}] \\
\boldsymbol{r}(t) &= \langle \cos t, \sin t, t \rangle, & &  \hspace*{-1cm} 0 \leq t \leq 5\pi  & \qquad &[\text{helix}].
\end{align*}
\end{Example}
Example~\ref{exparplot} gives illustrations of parametric plots in both 2D plane and 3D space. The \emph{wxMaxima} input commands and their corresponding outputs are presented in Table~\ref{table2}. A 2D plot of the polar curve cardioid $r(t) = 1 - \sin t$, $0 \leq t \leq 2\pi$, and a 3D plot of helix $\boldsymbol{r}(t) = \langle \cos t, \sin t, t \rangle$, $0 \leq t \leq 5\pi$ are given on the left and right panels of Figure~\ref{fig3-para}, respectively. 
\end{paracol}
\nointerlineskip
\begin{paracol}{2}
\switchcolumn

Two of integral applications are calculating the area enclosed by a polar curve and the arc length of a parametric curve. The commands to calculate the area $A$ enclosed by the cardioid and the length of the helix $L$ and their corresponding results are also given in Table~\ref{table2}.\linebreak Mathematically, they are performed as follows, where the prime denotes differentiation with respect to $t$:
\begin{align*}
A &= \int_{0}^{2 \pi} \frac{1}{2} r^2 \,dt = \int_{0}^{2 \pi} \frac{1}{2} (1 - \sin t)^2 \, dt = \frac{3\pi}{2} \\
L &= \int_{0}^{5 \pi} \sqrt{x'(t)^2 + y'(t)^2 + z'(t)^2} \,dt = \int_{0}^{5 \pi} \sqrt{\sin^2 t + \cos^2 t + 1} \, dt = 5\sqrt{2}\pi.
\end{align*}

\begin{specialtable}[H]
\renewcommand{\arraystretch}{1.5}
\caption{Examples of \emph{wxMaxima} commands and their corresponding outputs showing 2D and 3D plots of parametric curves cardioid and helix, respectively, and their corresponding arc lengths.}	\label{table2}
\begin{tabular*}{0.75\textwidth}{@{\extracolsep{\fill}} ll @{}}
\toprule
\textbf{Input }& \textbf{Output }														\\ 
\midrule
\verb|(%i10) r: 1 - sin(t);|									& ($r$) $1 - \sin(t)$	\\
\verb|(%i11) plot2d([parametric,r*cos(t),r*sin(t)],| 			& (Figure~\ref{fig3-para},  \\
\verb|       [t,0,2*%pi],[color,red],[x,-1.5,1.5],|  			& left panel) \\
\verb|       [y,-2.25,0.5],same_xy,[title,"A cardioid"]);| 		& \\
\verb|(%i12) integrate(1/2*r^2,t,0,2*%pi);|  					& ${\displaystyle \frac{3\pi}{2}}$ \\
\verb|(%i13) plot3d([cos(t),sin(t),t],[t,0,5*%pi],[y,-1,1],| 	&  \\
\verb|       [grid,100,2],[gnuplot_pm3d,true],[elevation,50],| 	& (Figure~\ref{fig3-para},  \\
\verb|       [azimuth,60],[legend,false],[title,"A helix"]);| 	& right panel) \\
\verb|(%i14) factor(integrate(sqrt(diff(cos(t),t)^2| 			& \\
\verb|       +diff(sin(t),t)^2+diff(t,t)^2),t,0,5*%pi));| 		& $5 \sqrt{2} \pi$ \\
\bottomrule
\end{tabular*}

\end{specialtable}
\begin{Example}
3D surfaces of M\"{o}bius band and torus:
\begin{align*}
\boldsymbol{s}(x,y) &= \left\langle \left(3 + y \cos \frac{x}{2} \right) \cos x, \left(3 + y \cos \frac{x}{2}\right) \sin x, y \sin \frac{x}{2} \right\rangle & & [\text{M\"{o}bius band}]\\
\boldsymbol{t}(\theta,\phi) &= \langle (2 + \cos \theta) \cos \phi, (2 + \cos \theta) \sin \phi, \sin \theta \rangle, \qquad 0 \leq \theta, \phi < 2\pi. & \qquad & [\text{torus}]
\end{align*}	
\end{Example}

The \emph{wxMaxima} input commands and their corresponding outputs are given in Table~\ref{table3} and Figure~\ref{fig4-tor}, respectively.  Three-dimensional plots such as this M\"{o}bius band and torus can be rotated easily in any direction. Sometimes also called the M\"obius strip, the former is a surface with only one side and only one boundary curve, the simplest example of a nonorientable surface. All M\"obius bands have a twofold symmetry rotational axis, for which a 180$^{\circ}$ rotation results in strips indistinguishable from the original. The surface finds abundant applications in physical sciences, including nanostructures~\cite{straosin}, metamaterials~\cite{chang}, polymeric materials~\cite{nie}, and photonic crystals~\cite{han}. 

The torus example comes from the problem of calculating the volume of solid revolution. It is obtained from a disk $(x - 2)^2 + z^2 \leq 1$ revolved about the $z$-axis, and thus, it is radially symmetric about the $z$-axis. The desired viewpoint can be obtained easily by setting the elevation and azimuth angles. This is beneficial in comparison to the \emph{Matlab} plots where figures usually have a high number of pixels and are rather heavy to be rotated. Although the surface of revolution is a fascinating object among geometers and topologists, it also finds ample applications in nanophotonics~\cite{ahmadivand}, metamaterials~\cite{kaelberer}, magnetized plasma~\cite{kliem}, and polymer chemistry~\cite{pochan}.
\begin{specialtable}[H]
\renewcommand{\arraystretch}{1.25}
\caption{Examples of \emph{wxMaxima} commands and their corresponding outputs showing 3D plots of parametric surfaces M\"obius band and torus.}	\label{table3}
\begin{tabular*}{0.75\textwidth}{@{\extracolsep{\fill}} ll @{}}
\toprule
\textbf{Input} 	& \textbf{Output}											  \\ 
\midrule
\verb|(%i15) plot3d([cos(x)*(3+y*cos(x/2)),sin(x)*(3+y*cos(x/2)),|	& \\
\verb|       y*sin(x/2)],[x,-%pi,%pi],[y,-1,1],['grid,50,15],| 		& (Figure~\ref{fig4-tor},  \\
\verb|       [legend,false],[elevation,35],[azimuth,50],|  			& left panel) \\
\verb|       [title,"Moebius band"]);| 								& \\
\verb|(%i16) plot3d([cos(y)*(2+cos(x)),sin(y)*(2+cos(x)),sin(x)],|  & \\
\verb|       [x,0,2*%pi],[y,0,2*%pi],[gnuplot_pm3d,true],| 			& (Figure~\ref{fig4-tor},\\
\verb|       [grid,50,50],[legend,false],[elevation,30],| 			& right panel)  \\
\verb|       [azimuth,135],[title,"A torus"]);| 					& \\
\bottomrule
\end{tabular*}
\end{specialtable}

\end{paracol}
\nointerlineskip
\unskip
\begin{paracol}{2}
\switchcolumn

\subsection{\emph{wxMaxima}'s Weaknesses}

In this subsection, we consider several examples where \emph{wxMaxima} is deficient. When sketching a function with singularities, we can discern an improved plot by including supplementary commands. The CAS \emph{wxMaxima} is not particularly strong when it comes to evaluating integrals, particularly the integral of rational functions. This may come as a surprise since the latter are commonplace and considered to be schematically viable.
\begin{Example}
A sketch of ${\displaystyle y = f(x) = \frac{x^2}{x^2 - 1}}$.
\end{Example}	
In this example, sketching a 2D plot of a function with vertical asymptotes requires an additional instruction. Although the output from ``\verb|plot2d(x^2/(x^2-1),[x,-4,4]);|'' input shows that there are two vertical asymptotes at $x = \pm 1$, modifying the input produces a better visibility,~cf.~\cite{glasner}:\\ ``\verb|plot2d(max(min(x^2/(x^2-1),10),-10),[x,-4,4]);|''. Figure~\ref{fig4} compares the 2D plots obtained from these two inputs.

\end{paracol}
\nointerlineskip
\vspace{-12pt}
\begin{paracol}{2}
\switchcolumn

\begin{Example}		\label{extrig}
Indefinite and definite integrals involving trigonometric functions:
\begin{equation*}
\int \frac{x \sin x}{1 + \cos^2 x} \, dx \qquad \qquad \text{and} \qquad \qquad \int_{0}^{\pi} \frac{x \sin x}{1 + \cos^2 x} \, dx.
\end{equation*}
\end{Example}
Both \emph{wxMaxima} and \emph{Matlab} fail to evaluate both integrals. The indefinite integral above involves a non-elementary function Li$_2(z)$.  Special function Li$_s(z)$ is a polylogarithm of order $s$ and argument $z$. It is also known as the Jonqui\`ere's function and is defined by a power series in~$z$; it is also a Dirichlet series in~$s$:
\begin{equation*}
\text{Li}_s(z) = \sum_{k = 1}^{\infty} \frac{z^k}{k^s}.
\end{equation*}
Both \emph{WolframAlpha} and \emph{Mathematica} can produce a correct answer for the above definite integral: $\pi^2/4$. Although treatment of the indefinite integral is beyond the Calculus course, the definite integral can be evaluated using a simple substitution rule. Appendix~\ref{appeA} shows how to evaluate it without any help from~CAS.

\begin{Example}
G\"artner's indefinite integral of a rational function~\cite{gaertner}:
\begin{equation*}
\int \frac{\sqrt{2}\, dx}{(x - 1)^4 + \frac{1}{16}}.
\end{equation*}
\end{Example}
When employing the  ``\verb|integrate(sqrt(2)/((x-1)^4+1/16),x);|'' command, \emph{wxMaxima} could not directly handle this integral. However, using the substitution $y = x - 1$:
``\verb|changevar(%,x-1-y,y,x);|'', the integral transforms into
\begin{equation*}
\int \frac{16\sqrt{2}\, dy}{16y^4 + 1}.
\end{equation*}
Asking \emph{wxMaxima} to evaluate the integral, ``\verb|ev(%,integrate);|'', we obtain
\begin{equation*}
2 \log(4y^2 + 2 \sqrt{2}y + 1) - 2 \log(4y^2 - 2 \sqrt{2}y + 1) + 4 \atan \left(\frac{4y + \sqrt{2}}{\sqrt{2}} \right) + 4 \atan \left(\frac{4y - \sqrt{2}}{\sqrt{2}} \right).
\end{equation*}
\textls[-15]{Substituting back to the original variable gives us the desired integral,} ``\verb|sfx: %,y=x-1|'':
\begin{align*}
& 2 \log(4(x - 1)^2 + 2 \sqrt{2}(x - 1) + 1) - 2 \log(4(x - 1)^2 - 2 \sqrt{2}(x - 1) + 1) \\
& + 4 \atan \left(\frac{4(x - 1) + \sqrt{2}}{\sqrt{2}} \right) + 4 \atan \left(\frac{4(x - 1) - \sqrt{2}}{\sqrt{2}} \right).
\end{align*}

In attempting to solve indefinite integrals of rational functions, \emph{wxMaxima} essentially still uses the algorithms that were described by Moses~\cite{moses}. The CAS employs the classical algorithmic tool in symbolic integration known as Hermite reduction. Both Hermite and Ostrogradsky demonstrated decomposing indefinite integral of rational functions as the sum of a function with simple poles and the derivative of another rational function, i.e., for $R \in \mathbb{Q}(x)$
\begin{equation*}
\int R(x) \, dx = U(x) + \int A(x) \, dx, \qquad  U, \; A \; \in \; \mathbb{Q}(x)
\end{equation*}
where $A$ has only simple poles and vanishes at infinity~\cite{subramaniam,bostan13,bostan18,moir}.

Thus, \emph{wxMaxima} implements the combination of Hermite reduction with either the partial fraction decomposition (Bernoulli--Leibniz algorithm) to compute the rational part of the integrand or with some heuristic methods in computing its transcendental parts for several easier cases. Although other commercial software such as \emph{Axiom} and \emph{3M} might also implement integration using Hermite reduction for the rational part of the integrand, it employs more contemporary and advanced methods to integrate the remaining terms after Hermite reduction is implemented.
	
\begin{Example}	\label{exbrons}
Bronstein's definite integral of a rational function~\cite{bronstein}:
\begin{equation*}
\int_{1}^{2} \frac{x^4 - 3x^2 + 6}{x^6 - 5 x^4 + 5 x^2 + 4} \, dx.
\end{equation*}
\end{Example}
When we attempt the command ``\verb|integrate((x^4-3*x^2+6)/(x^6-5*x^4+5*x^2+4),|\\ \verb|x,1,2);|'', \emph{wxMaxima} could not compute this integral. However, \emph{Axiom} and \emph{3M} can tackle it. \emph{Axiom}'s output reads:
\begin{equation*}
\frac{2 \atan(8) + 2 \atan(5) + 2 \atan(2) + 2 \atan\left(\frac{1}{2}\right) - \pi}{2}.
\end{equation*}
\emph{Maple}'s output is given as follows:
\begin{equation*}
\arctan\left(\frac{1}{2}\right) - \frac{1}{2} \pi + \arctan(5) + \arctan(8) + \arctan(2).
\end{equation*}
Interestingly, \emph{Mathematica} and \emph{Matlab} give the following outputs, respectively, which were simpler than the ones from the previous two software:
\begin{equation*}
\frac{5}{4} \pi - \arctan(2) \qquad \text{and} \qquad  \pi - \arctan(1/3).
\end{equation*} 
All of these results correspond to an identical numerical value of approximately $2.819842$. The proof is shown in Appendix~\ref{appeB1}.

It would also be interesting to observe the corresponding indefinite integral. While \emph{wxMaxima} anticipatedly cannot compute the indefinite integral,  \emph{Axiom}, \emph{Maple}, and \emph{Matlab} delivered the same results and agreed with Formula (2.20) in~\cite{bronstein}. This reads
\begin{equation*}
\int \frac{x^4 - 3x^2 + 6}{x^6 - 5 x^4 + 5 x^2 + 4} \, dx = \arctan (x) + \arctan (x^3) + \arctan \left(\frac{x^5 - 3x^3 + x}{2} \right) + C.
\end{equation*}
In particular, \emph{Mathematica} is of special interest as the output of the corresponding indefinite integral does not agree with the formula given in~\cite{bronstein}. It reads as follows (some terms were rearranged):
\begin{equation*}
\int \frac{x^4 - 3x^2 + 6}{x^6 - 5 x^4 + 5 x^2 + 4} \, dx = \frac{1}{2} \arctan \left[\frac{x(x^2 -3)}{x^2 - 2}\right] - \frac{1}{2} \arctan \left[\frac{x(x^2 - 3)}{2 - x^2} \right] + C.
\end{equation*}
\noindent
The definite integral is given as follows:

\end{paracol}
\nointerlineskip
\begin{align*}
\int_{1}^{2} \frac{x^4 - 3x^2 + 6}{x^6 - 5 x^4 + 5 x^2 + 4} \, dx &= 
\lim\limits_{t \to \sqrt{2}^{{\,}-}} \left \{ \frac{1}{2} \arctan \left[\frac{x(x^2 -3)}{x^2 - 2}\right] - \frac{1}{2} \arctan \left[\frac{x(x^2 - 3)}{2 - x^2} \right] \right \} - \arctan (2) \\
& \qquad + \arctan (1) - \lim\limits_{t \to \sqrt{2}^{+}} \left \{ \frac{1}{2} \arctan \left[\frac{x(x^2 -3)}{x^2 - 2}\right] - \frac{1}{2} \arctan \left[\frac{x(x^2 - 3)}{2 - x^2} \right] \right \} \\
&= \left(\frac{\pi}{2} - \arctan (2) \right) + \left(\frac{\pi}{4} + \frac{\pi}{2} \right) = \frac{5}{4} \pi - \arctan (2).
\end{align*}
\begin{paracol}{2}
\switchcolumn
The first and second chapters of~\cite{bronstein} explain outline algorithms for the integration of rational functions in great detail. All relevant algorithms are given in pseudocodes. With an estimated effort of 500 to 700 lines of code, these algorithms can be implemented in \emph{wxMaxima}.

\begin{Example}	\label{exadam}
Adamchik's definite integral of a rational function~\cite{adamchik}:
\begin{equation*}
\int_{0}^{4} \frac{x^2 + 2x + 4}{x^4 - 7x^2 + 2x + 17} \, dx.
\end{equation*}
\end{Example}
For an indefinite integral, \emph{wxMaxima} could not compute it. Other software, namely, \emph{Axiom}, \emph{Maple}, and \emph{Matlab}, presented an identical solution but with rather dissimilar expressions: 
\begin{equation*}
\int \frac{x^2 + 2x + 4}{x^4 - 7x^2 + 2x + 17} \, dx = \arctan(x - 1) + \arctan\left(\frac{1}{3} x^3 - \frac{1}{2} x^2 - x + \frac{5}{3} \right) + C.
\end{equation*}
On the other hand, the computational result from \emph{Mathematica} produced an antiderivative expression with discontinuities at $x = \pm 2$:
\begin{equation*}
\int \frac{x^2 + 2x + 4}{x^4 - 7x^2 + 2x + 17} \, dx = \frac{1}{2} \arctan \left(\frac{-x - 1}{x^2 - 4} \right) - \frac{1}{2} \arctan\left(\frac{x + 1}{x^2 - 4} \right) + C.
\end{equation*}
For the definite integral, while all numerical results agree, the symbolic outputs yield remarkably distinct expressions. The outputs are given as follows:
\end{paracol}
\clearpage 
\nointerlineskip
\begin{align*}
\int_{0}^{4} \frac{x^2 + 2x + 4}{x^4 - 7x^2 + 2x + 17} \, dx &= \frac{\pi}{4} + \arctan (3) - \arctan \left(\frac{5}{3}\right) + \arctan \left(\frac{41}{3}\right)  & & \text{[\emph{Axiom} and \emph{Maple}]} \\
&= \frac{5}{4}\pi - \arctan \left(\frac{75}{11} \right) & & \text{[\emph{Matlab}]} \\
&= \pi - \arctan \left(\frac{1}{4} \right) - \arctan \left(\frac{5}{12}\right) & &\text{[\emph{Mathematica}]} \\
&= (\text{a tedious and complicated expression})  & &\text{[\emph{wxMaxima}]} \\
&\approx 2.50182  & & \text{[numerical value]}.
\end{align*}
\begin{paracol}{2}
\switchcolumn

Although \emph{wxMaxima} computes and produces an output for the definite integral, that expression is exceptionally tedious and complicated. It is also unclear how it was obtained, but possibly by performing a contour integration. Despite that tiresome and intricate expression, its numerical evaluation yields a value that agrees with the one calculated by other software. A high floating-point precision command, e.g., ``\verb|fpprec:150|'' needs to be imposed to obtain that numerical value. Except for the exact \emph{wxMaxima} output, Appendix~\ref{appeB2} verifies that the above outputs had identical values.

Similar to Example~\ref{exbrons}, when we want to use \emph{Mathematica}'s antiderivative result to compute the value of the definite integral, we need to split the interval of integration at the point of discontinuity to avoid a wrong result of ${\displaystyle -\arctan \left(\frac{1}{4}\right) - \arctan \left(\frac{5}{12}\right)}$. The computational summary is given as follows:
\end{paracol}
\nointerlineskip
\begin{align*}
\int_{0}^{4} \frac{x^2 + 2x + 4}{x^4 - 7x^2 + 2x + 17} \, dx &= 
\lim\limits_{t \to 2^{{\,}-}} \left \{ \frac{1}{2} \arctan \left(\frac{-t - 1}{t^2 - 4} \right) - \frac{1}{2} \arctan\left(\frac{t + 1}{t^2 - 4} \right) \right \} - \arctan \left(\frac{1}{4} \right) \\
& \qquad - \arctan \left(\frac{5}{12} \right) - \lim\limits_{t \to 2^{+}} \left \{ \frac{1}{2} \arctan \left(\frac{-t - 1}{t^2 - 4} \right) - \frac{1}{2} \arctan\left(\frac{t + 1}{t^2 - 4} \right) \right \} \\
&= \left[\frac{\pi}{2} - \arctan \left(\frac{1}{4} \right) \right] - \left[\arctan \left(\frac{5}{12} \right) - \frac{\pi}{2} \right] \\
&= \pi - \arctan \left(\frac{1}{4} \right) - \arctan \left(\frac{5}{12} \right).
\end{align*}
\begin{paracol}{2}
\switchcolumn

\begin{Example}
Tobey's indefinite integral of a rational function~\cite{tobey}:
\begin{equation*}
\int \frac{7x^{13} + 10 x^8 + 4x^7 - 7x^6 - 4x^3 - 4x^2 + 3x + 3}{x^{14} - 2 x^8 - 2x^7 - 2x^4 - 4x^3 - x^2 + 2x + 1} \, dx.
\end{equation*}
\end{Example}\vspace{12pt}
Except for \emph{wxMaxima}, all other CAS can compute this integral. The result can be found on page~502 of~\cite{geddes} using the Rothstein--Trager method. It reads
\end{paracol}
\begin{equation*}
\frac{1}{2}(1 + \sqrt{2}) \log \left(x^7 - \sqrt{2} \, x^2 - (1 + \sqrt{2}) x - 1\right) +
\frac{1}{2}(1 - \sqrt{2}) \log \left(x^7 + \sqrt{2} \, x^2 - (1 - \sqrt{2}) x - 1\right) + C.
\end{equation*}
\begin{paracol}{2}
\switchcolumn

\clearpage
\section{Discussion\\} \label{discussion} 

\subsection{\emph{wxMaxima} in Perspective}

There exists a distinct paradigm between the CAS mentioned in the Introduction. While \emph{Matlab} is particularly strong in numerical computation, the real forte of \emph{wxMaxima}, \emph{Maple}, and \emph{Mathematica} is symbolic computation. \emph{Matlab} has also incorporated some symbolic features, and the latter three can also perform numerical computations as well.

When it comes to symbolic functionality, \emph{wxMaxima} still has a long way to go to catch up  other commercial software. Together with \emph{Axiom} and \emph{SageMath}, a formula editor is missing in \emph{Maxima}, although \emph{wxMaxima} serves efficiently as its user interface. While \emph{Axiom} does not possess graph-theory symbolic functionality, \emph{wxMaxima} does. Unfortunately, Diophantine equation solvers and quantifier elimination are absent in \emph{wxMaxima}, and \emph{SageMath} acquires these two features via \emph{SymPy} and \emph{qepcad} optional packages, respectively. Table~\ref{functionality} summarizes the significantly developed symbolic functionality in various software systems.
\end{paracol}
\nointerlineskip
\begin{specialtable}[H]
\widetable
\caption{A summary of significantly developed symbolic functionality among various software. \emph{Maxima} has a shortfall in quantifier elimination and Diophantine equation solver. The missing formula editor can be overcome by a user interface \emph{wxMaxima}.}	\label{functionality}}
{\footnotesize
\tabcolsep=0.14cm 
\begin{tabular}{@{}lcccccccc@{}}
\toprule
		 			& \multirow{3}{*}{{\textbf{Formula Editor}}}  & \multicolumn{2}{c}{\textbf{Calculus}} 	& \multirow{3}{*}{{\textbf{Quantifier Elimination}}}  & \multicolumn{4}{c}{\textbf{Solvers}} 	 		 	\\ 
\cline{3-4} \cline{6-9}
\textbf{Software	}		& 		   & \multirow{2}{*}{\textbf{Integration}} & \textbf{Integral 	  }	& 	& \textbf{Inequalities}& \textbf{Diophantine} 	& \textbf{Differential }	& \textbf{Recurrence}	  \\ 
					&		   &		     & \textbf{Transforms }	  	& 				&    			& \textbf{Equations}			& \textbf{Equations} 	& \textbf{Relations }				\\ 
\midrule
\emph{Axiom} 		& \xmark   & \cmark		& \cmark 		  	& \cmark  		& \cmark 		& \cmark & \cmark		& \cmark			\\
\emph{Magma}		& \xmark   & \xmark 	& \xmark 		  	& \xmark		& \cmark		& \xmark & \xmark		& \xmark			\\
\emph{Maple}		& \cmark   & \cmark		& \cmark 		  	& \cmark 		& \cmark		& \cmark & \cmark		& \cmark			\\
\emph{Mathematica}	& \cmark   & \cmark 	& \cmark 			& \cmark 		& \cmark		& \cmark & \cmark		& \cmark			\\
\emph{Maxima}		& \xmark   & \cmark		& \cmark 			& \xmark 		& \cmark		& \xmark & \cmark		& \cmark			\\
\emph{Matlab}		& \cmark   & \cmark		& \cmark 			& \cmark 		& \cmark		& \xmark & \cmark		& \xmark 			\\
\emph{SageMath}		& \xmark   & \cmark		& \cmark 			& \cmark 		& \cmark		& \cmark & \cmark		& \cmark			\\
\bottomrule
\end{tabular}
\end{specialtable}
\begin{paracol}{2}
\switchcolumn

An important feature of \emph{wxMaxima} that trumps other CAS is the admittance of referring to the result of the last evaluated expression with a percentage sign~(\verb|%|). Although we may organize the commands in spatial order, \emph{wxMaxima} stores information in chronological order. Thus, the~\verb|%| always refers to the most recently executed command, and not necessarily to the one that appears directly above the executed command~\cite{timberlake}. As already mentioned in the Introduction, since \emph{wxMaxima} is open-source software, it is maintained by an active community of developers, and thus being updated regularly. New releases occur approximately twice annually. For a proprietary software like \emph{3M}, not everyone can inspect, modify, or enhance, but only a team of developers from the company or organization who maintains exclusive control over it can modify the software. Admittedly, the most important reason for utilizing \emph{wxMaxima}, both in teaching and research is the cost, and many authors agreed on this particular notable reason~\cite{hannan,timberlake,senese}.

\subsection{Educational Benefit}	

There are many educational benefits of embedding CAS \emph{wxMaxima} into teaching and learning, Calculus in particular, and mathematics in general. The software serves as meaningful assistance not only to students but also to instructors in verifying hand-computed results, obtaining numerical values of  tedious expressions, and providing rich-quality graphical plots. Hence, the benefit of algebraic, numerical, and graphical aspects, respectively. The time reduction gained in the labor of manual calculation can be spent for a deeper understanding of key ideas, theoretical concepts, and problem-solving methods.

Some examples discussed in Section~\ref{strongweak} support this testimony. By observing the graph of a cubic function, students should establish a connection between the visual plot and the computational result of obtaining intervals where the function is increasing, decreasing, concave up, and concave down. Other features, such as the local maxima and minima, and inflection points should follow  naturally. By examining the linear approximation of a function, students should be able to conclude why the approximation works for the variable $x$ near a specified point of interest. At an advanced level, this comprehension is advantageous for linear-stability analysis. 

The literature does not lack in supplying educational benefits of using CAS in mathematics lessons. From a coding perspective, the programming language in \emph{wxMaxima} flows more naturally in comparison to that in other software, and hence allows for students to autonomously implement simple algorithms~\cite{garcia}. For current and future calculus teachers, \emph{wxMaxima} could reduce the time spent on course preparation~\cite{velychko}. By allowing students to induce judgments and make mistakes, the adoption of \emph{wxMaxima} stimulates the interactive-learning process through testing, evaluation, decision-making, and error correction~\cite{zakova}. Weigand argued that the wise use of CAS could foster students' ability in problem solving, modeling, proving, and communicating~\cite{weigand}. The CAS can even be blended with an innovative pedagogical approach, such as flipped classrooms~\cite{karjanto19}.

Additional benefits come along as more instructors adopt and implement a CAS, particularly \emph{wxMaxima}, into their teaching and learning. However,  \emph{wxMaxima} is not a perfect software, as we elaborated through several examples in Section~\ref{strongweak}. 

\subsection{Limitation} \label{limitation}

This study admits several limitations. First, the extent of course materials in Single-Variable Calculus (SVC) is overwhelmingly profuse while the time is scanty. It covers nine chapters of Stewart's Calculus textbook~\cite{stewart}, and they need to be completed in 14 weeks. Embedding technology into teaching and learning contributes an additional burden. In a typical North American university, similar content would be covered in two different courses for two consecutive semesters. Usually called ``Calculus~1'' and ``Calculus~2'', they cover Differential and Integral Calculus, respectively. The latter often includes an Introduction to Differential Equations, and Sequence and Series. Hence, teaching Calculus using \emph{wxMaxima} with less material coverage seems to be a promising attempt, as we could pursue in MVC.

Second, despite the superior features of \emph{wxMaxima} and its specialty in symbolic operations, the CAS itself has some weaknesses that could be challenging for beginners to learn, adopt, and adapt, both as instructors and students, where it might be more conspicuous for the latter. The CAS is admittedly far from perfect due to the nature of free and open-source software. Other commercial CAS such as \emph{Mathematica} and \emph{Maple} might be better since the companies that release them possess an army of paid personnel working around the clock to improve the software. Nevertheless, we are not without hope since \emph{wxMaxima} is updated frequently, bugs are fixed, and documentation is improved by a group of volunteer developers who work tirelessly.

\section{Conclusions}	\label{conclusion}

In this article, we considered several features of \emph{wxMaxima} that could be useful for enhancing the quality of Calculus teaching and learning. Although the CAS itself is far from perfect, focusing on its strengths might benefit both students and instructors when embedding the software into the subject. We included some examples where \emph{wxMaxima} assists well in understanding Calculus concepts better. These include calculating limit, finding a derivative of a function, evaluating definite and indefinite integrals, generating plots for explicit functions, parametric functions, polar curves, and three-dimensional objects. The visualization aspect enhances excellent teaching.

For teachers and instructors who are currently adopting \emph{Maple} or \emph{Mathematica} in their teaching, we are interested in stimulating a discussion on whether it is viable to integrate, or even switch entirely, to \emph{wxMaxima}. This endeavor is not beyond our reach, at least for a sequence of Calculus courses (PreCalculus, SVC, and MVC).

For Linear Algebra, \emph{Matlab} is still the most popular software among both engineers and educators thanks to its many professional contributors, rigorous development, powerful numerical computation, and additional package \emph{Simulink}. \emph{SageMath} is also becoming popular for Linear Algebra teaching and learning since it utilizes existing open-source libraries specifically designed for Linear Algebra, including \emph{LAPack} and \emph{NumPy}. Since \emph{wxMaxima} also possesses many functions for manipulating matrices, we would not be too ambitious to persuade educational practitioners by considering to switch CAS to \emph{wxMaxima} for Linear Algebra teaching and learning. Some Linear Algebra problems are worth testing using \emph{wxMaxima} nonetheless.

Despite the many admirable qualities of \emph{wxMaxima}, our experience in other courses and with different CAS suggests that some students hate when  instructors generally attempt to embed any CAS in general into mathematics teaching and learning. For example, after implementing Calculus with \emph{SageMath} in 2012, some students provided feedback suggesting to eliminate \emph{SageMath} from Calculus teaching. For Linear Algebra, a dedicated one-hour problem-solving session using \emph{Matlab} was not favored among the students. Students' feedback on teaching evaluation consistently mentions that the computer laboratory sessions are a waste of time and must be replaced by the traditional problem-solving sessions using pen and paper.

From our experience of embedding \emph{wxMaxima} into Calculus courses, both Single-Variable and Multivariable, very few students gave positive feedback regarding the software. The majority of students' comments voice an atmosphere of negativity and resistance, they tend to push away \emph{wxMaxima}. Even if they seem to embrace the CAS favorably, they might forget it as soon as the semester is over. In subsequent courses, their instructors might not use \emph{wxMaxima} anymore. There are other important programming languages (e.g., \emph{Python}, \emph{C++}, and \emph{Java}) that students could master before they enter the job market.

\vspace{6pt}

\authorcontributions{Conceptualization, N.K.; methodology, N.K.; software, H.S.H.; validation, N.K. and H.S.H.; formal analysis, N.K.; investigation, N.K.; resources, N.K.; data curation, N.K.; writing---original draft preparation, N.K.; writing---review and editing, N.K.; visualization, N.K.; supervision, N.K.; project administration, N.K.; funding acquisition, N.K. All authors have read and agreed to the published version of the manuscript.}

\funding{This research received no external funding.}

\institutionalreview{Not applicable.}

\informedconsent{Not applicable.}

\dataavailability{Not applicable.}

\acknowledgments{The authors gratefully acknowledge Boris G\"artner from Munich, Germany for the fruitful discussion on the limitations of \emph{wxMaxima} and for providing references to the relevant literature, particularly~\cite{moses} and~\cite{bronstein,adamchik,tobey,geddes}.}

\conflictsofinterest{The authors declare no conflict of interest.}

\section*{\small{Dedication}}
The main author would like to dedicate this article to the memory of his late father Zakaria~Karjanto (Khouw~Kim~Soey, 許金瑞) who not only taught him the alphabet, numbers, and the calendar in his early childhood, but also cultivated the value of hard work, diligence, discipline, perseverance, persistence, and grit. Karjanto Senior was born in Tasikmalaya, West~Java, Japanese-occupied Dutch~East~Indies on 1~January~1944 (Saturday~Pahing) and died in Bandung, West~Java, Indonesia on 18~April~2021 (Sunday~Wage).

\appendixtitles{yes}
\appendixstart
\appendix

\section{Evaluating a Definite Integral without CAS}	\label{appeA}
The definite integral involving trigonometric functions considered in Example~\ref{extrig} can be evaluated without the aid of any CAS. We have the following lemma:
\begin{Lemma}	\label{lemma1}
For a continuous and rational function $f(x) \in \mathbb{Q}(x)$
\begin{equation*}
\int_{0}^{\pi} x \, f(\sin x) \, dx = \frac{\pi}{2} \int_{0}^{\pi} f(\sin x) \, dx.
\end{equation*}	
\end{Lemma}
\begin{proof}
Let $u = \pi - x$, then $du = -dx$, for $x = 0$, $u = \pi$ and $x = \pi$, $u = 0$. We have
\begingroup
\makeatletter\def\f@size{9}\check@mathfonts
\def\maketag@@@#1{\hbox{\m@th\normalsize\normalfont#1}}%
\begin{align*}
	\int_{0}^{\pi} x \, f(\sin x) \, dx &= \int_{\pi}^{0} (\pi - u) \, f(\sin (\pi - u)) \, (-du) 
	= -\pi \int_{\pi}^{0} f(\sin u) \, du + \int_{\pi}^{0} u \, f(\sin u) \, du \\
	&=  \pi \int_{0}^{\pi} f(\sin x) \, dx - \int_{0}^{\pi} x \, f(\sin x) \, dx \\
	2 \int_{0}^{\pi} x \, f(\sin x) \, dx &= \pi \int_{0}^{\pi} f(\sin x) \, dx \\
	\int_{0}^{\pi} x \, f(\sin x) \, dx &= \frac{\pi}{2} \int_{0}^{\pi} f(\sin x) \, dx.
	\end{align*}
\endgroup
The proof is complete.
\end{proof}	
We now have the following proposition.
\begin{Proposition}
\begin{equation*}
\int_{0}^{\pi} \frac{x \sin x}{1 + \cos^2 x} \, dx = \frac{\pi^2}{4}.
\end{equation*}	
\end{Proposition}
\begin{proof}
There are at least three approaches with subtle differences in tackling this problem. The first method is by taking $f(x) = x/(2 - x^2)$ and applying Lemma~\ref{lemma1}. We then obtain
\begingroup\makeatletter\def\f@size{8}\check@mathfonts
\def\maketag@@@#1{\hbox{\m@th\normalsize\normalfont#1}}%
\begin{align*}
\int_{0}^{\pi} x f(\sin x) \, dx &= \int_{0}^{\pi} \frac{x \sin x}{2 - \sin^2 x} \, dx 
 = \int_{0}^{\pi} \frac{x \sin x}{1 + \cos^2 x} \, dx 
 = \frac{\pi}{2} \int_{0}^{\pi} \frac{\sin x}{1 + \cos^2 x} \, dx \\
&=-\frac{\pi}{2} \int_{0}^{\pi} \frac{d(\cos x)}{1 + \cos^2 x} \, dx
 = \left. -\frac{\pi}{2} \tan^{-1} (\cos x) \right|_{0}^{\pi}
= -\frac{\pi}{2} \left(\tan^{-1} (-1) - \tan^{-1} 1 \right) = \frac{\pi^2}{4}.
\end{align*}
\endgroup
The second technique is by writing $x = \left(x - \frac{\pi}{2}\right) + \frac{\pi}{2}$ and substituting $u = x - \frac{\pi}{2}$. It becomes
\begin{equation*}
\int_{0}^{\pi} \frac{x \sin x}{1 + \cos^2 x} \, dx = 
\int_{-\frac{\pi}{2}}^{\frac{\pi}{2}} \frac{u \cos u}{1 + \sin^2 u} \, du + \frac{\pi}{2} \int_{0}^{\pi} \frac{\sin x}{1 + \cos^2 x} \, dx.
\end{equation*}
The first integral of the right-hand side vanishes since the integrand is an odd function. The second integral follows the first method. The third approach is by substituting $u = x - \frac{\pi}{2}$ directly from the beginning. It yields a slightly different expression from the previous two approaches for the second term of the right-hand side:
\begin{equation*}
\int_{0}^{\pi} \frac{x \sin x}{1 + \cos^2 x} \, dx = 
\int_{-\frac{\pi}{2}}^{\frac{\pi}{2}} \frac{u \cos u}{1 + \sin^2 u} \, du + \frac{\pi}{2} \int_{-\frac{\pi}{2}}^{\frac{\pi}{2}} \frac{\cos u}{1 + \sin^2 u} \, du.
\end{equation*}
Similar to the second technique, the first integral on the right-hand side is zero. Since the integrand of the second integral on the right-hand side is an even function, it simplifies to twice of the integral from $u = 0$ to $u = \frac{\pi}{2}$. Employing another substitution $y = \sin u$, we obtain the desired result:
\begin{equation*}
\int_{0}^{\pi} \frac{x \sin x}{1 + \cos^2 x} \, dx = 
\pi \int_{0}^{\frac{\pi}{2}} \frac{d(\sin u)}{1 + \sin^2 u} = \pi \tan^{-1} (\sin u) \Big|_{0}^{\frac{\pi}{2}} = \pi (\tan^{-1} - 0) = \frac{\pi^2}{4}.
\end{equation*}
The proof is completed.
\end{proof}

\section{CAS Output Comparison\\}		\label{appeB}

\subsection{Bronstein's Definite Integral of a Rational Function}	\label{appeB1}

Before verifying that all outputs in Example~\ref{exbrons} have the same value, we need the following identity.
\begin{Lemma}	\label{lemarctan}
For all $x \in \mathbb{R} \setminus\{0\}$
\begin{equation*}
\arctan x + \arctan \left(\frac{1}{x}\right) = \left\{
\begin{array}{rl}
 \frac{\pi}{2}, & \qquad \text{if} \quad x > 0 \\
-\frac{\pi}{2}, & \qquad \text{if} \quad x < 0.
\end{array}
\right.
\end{equation*}
\end{Lemma}
\begin{proof}
For all $x \in \mathbb{R} \setminus\{0\}$, let
\begin{equation*}
f(x) = \arctan x + \arctan \left(\frac{1}{x}\right).
\end{equation*}
Then, $f$ is differentiable for every $x \neq 0$ and 
\begin{equation*}
f'(x) = \frac{1}{1 + x^2} + \frac{\left(-\frac{1}{x^2}\right)}{1 + \frac{1}{x^2}} = 0.
\end{equation*}
Hence, $f$ is constant on each connected component of all $x \in \mathbb{R} \setminus\{0\}$. Since $f(1) = \dfrac{\pi}{4} + \dfrac{\pi}{4} = \dfrac{\pi}{2}$, we conclude that $f(x) = \dfrac{\pi}{2}$ for all $x > 0$. And since $f(-1) = -f(1) = -\dfrac{\pi}{2}$, it follows that $f(x) = -\dfrac{\pi}{2}$ for all $x < 0$. We have completed the proof.
\end{proof}
\begin{Proposition}
All outputs in Example~\ref{exbrons} are equivalent.
\end{Proposition}
\begin{proof}
\textls[-30]{To show that all outputs are equivalent, we use the identity ${\displaystyle \arctan (2) + \arctan \left(\frac{1}{2}\right) = \frac{\pi}{2}}$} from Lemma~\ref{lemarctan}. Hence, for \emph{Mathematica} output, we only need to verify that ${\arctan (8) +}$ \linebreak ${ \arctan (5) = \frac{5}{4} \pi -}$ $\arctan (2)$. The left-hand side can be calculated as follows:
\begin{equation*}
\arctan (8) + \arctan (5) = \arctan \left(\frac{8 + 5}{1 - 8 (5)}\right) =  \arctan \left(\frac{13}{-39}\right) = \arctan \left(-\frac{1}{3}\right) \text{mod} \pi.
\end{equation*}
Applying again a similar identity from Lemma~\ref{lemarctan} for $x = -3 < 0$, ${\displaystyle \arctan \left(-\frac{1}{3}\right) =}$\linebreak ${ \arctan (3) -\frac{\pi}{2}}$, we only need to show that
\begin{equation*}
\arctan(3) + \arctan(2) = \frac{5}{4} \pi + \frac{\pi}{2} = \frac{7}{4} \pi.
\end{equation*}
The left-hand side is calculated as previously:
\begin{equation*}
\arctan(3) + \arctan(2) = \arctan\left(\frac{3 + 2}{1 - 3(2)}\right) = \arctan \left(\frac{5}{-5}\right) = \arctan(-1) \text{mod} \pi = \frac{7}{4} \pi.
\end{equation*}
For \emph{Matlab} output, we need to show that ${\displaystyle \arctan (8) + \arctan (5) = \pi - \arctan \left(\frac{1}{3} \right)}$. Using the identity from Lemma~\ref{lemarctan} that ${\displaystyle \arctan (3) + \arctan \left(\frac{1}{3} \right) = \frac{\pi}{2}}$, the right-hand side becomes ${\displaystyle \frac{\pi}{2} + \arctan (3)}$. Combining  term $\arctan 3$ with one term on the left-hand side, we obtain either
\begin{align*}
\arctan (5) + \arctan (8) - \arctan (3) & = \arctan (5) + \arctan \left(\frac{8-3}{1 + 8(3)}\right) \\
& = \arctan (5) + \arctan \left(\frac{1}{5}\right) \text{mod} {\pi} = \frac{\pi}{2} \qquad \text{or} \\
    \arctan (8) + \arctan (5) - \arctan (3) & = \arctan (8) + \arctan \left(\frac{5-3}{1 + 5(3)}\right) \\ 
& = \arctan (8) + \arctan \left(\frac{1}{8}\right) \text{mod} {\pi} = \frac{\pi}{2}.
\end{align*}
Both identities are correct when we take zero remainder in the congruence relationship. This completes the proof.	
\end{proof}

\subsection{Adamchik's Definite Integral of a Rational Function}	\label{appeB2}

The results from Adamchik's definite integral of a rational function are verified by the following lemma.
\begin{Lemma}
All outputs in Example~\ref{exadam} are identical. 
\end{Lemma}
\begin{proof}
We show that all exact values are identical. In each case, we take the zero remainder whenever the congruence relationship appears. First, we verify that the \emph{Axiom}/\emph{Maple}'s and \emph{Matlab}'s results are identical, i.e., 
\begin{align*}
\frac{5}{4}\pi - \arctan \left(\frac{75}{11} \right) &= \frac{\pi}{4} + \arctan (3) - \arctan \left(\frac{5}{3}\right) + \arctan \left(\frac{41}{3}\right) \\
\pi - \arctan \left(\frac{75}{11} \right) &= \arctan \left(\frac{2}{9} \right) \text{mod}{\pi} + \arctan \left(\frac{41}{3}\right) \\
&= \arctan \left(\frac{2}{9} \right) + \frac{\pi}{2} - \arctan \left(\frac{3}{41} \right) \qquad \text{(by Lemma~\ref{lemarctan})}\\
&= \frac{\pi}{2} \arctan \left(\frac{2/9 - 3/41}{1 + (2/9) (3/41)}\right) \text{mod}{\pi} \\
\frac{\pi}{2} &= \arctan \left(\frac{75}{11} \right) + \arctan \left(\frac{11}{75} \right).
\end{align*}
Second, we show that the \emph{Matlab}'s and \emph{Mathematica} results are identical, i.e.,
\begin{equation*}
\frac{5}{4}\pi - \arctan \left(\frac{75}{11} \right) = \pi - \arctan \left(\frac{1}{4} \right) - \arctan \left(\frac{5}{12}\right).
\end{equation*}
Bringing $\pi$ to the left-hand side and gathering all inverse tangent terms to the right-hand side, we obtain
\begin{align*}
\frac{\pi}{4} &= \arctan \left(\frac{75}{11} \right) - \arctan \left(\frac{1/4 + 5/12}{1 - (1/4) (5/12)} \right) \text{mod}{\pi} \\
&= \arctan \left(\frac{75}{11} \right) - \arctan \left(\frac{32}{48} \right) \text{mod}{\pi} 
 = \arctan \left(\frac{75/11 - 32/43}{1 + (75/11) (32/43)} \right) \text{mod}{\pi} \\
&= \arctan 1.
\end{align*}
Finally, we confirm that the \emph{Axiom}/\emph{Maple}'s and \emph{Mathematica}'s results are identical, i.e.,
\begingroup\makeatletter\def\f@size{8.5}\check@mathfonts
\def\maketag@@@#1{\hbox{\m@th\normalsize\normalfont#1}}%
\begin{align*}
\pi - \arctan \left(\frac{1}{4} \right) - \arctan \left(\frac{5}{12}\right) &= \frac{\pi}{4} + \arctan (3) - \arctan \left(\frac{5}{3}\right) + \arctan \left(\frac{41}{3}\right) \\
\frac{3}{4} \pi - \arctan \left(\frac{32}{43}\right) \text{mod}{\pi} &= \arctan \left(\frac{2}{9} \right) \text{mod}{\pi} + \arctan \left(\frac{41}{3}\right) \\
\frac{3}{4} \pi &= \frac{\pi}{2} - \arctan \left(\frac{43}{32}\right) + \arctan \left(\frac{2}{9} \right) + \arctan \left(\frac{41}{3}\right) \qquad \text{(by Lemma~\ref{lemarctan})} \\
\frac{\pi}{4} &= \arctan \left(\frac{2}{9} \right) + \arctan \left(\frac{41/3 - 43/32}{1 + (41/3) (43/32)}\right) \text{mod}{\pi} \\
\frac{\pi}{4} &= \arctan \left(\frac{2}{9} \right) + \arctan \left(\frac{7}{11}\right) = \arctan 1.
\end{align*}
\endgroup
The proof is complete.
\end{proof}
\end{paracol}

\reftitle{References}

\end{document}